\documentclass[a4paper, 11pt]{article}

\usepackage{amsmath,amsfonts,latexsym,theorem}
\nonstopmode

\numberwithin{equation}{section}
\renewcommand{\theenumi}{\roman{enumi}}
\renewcommand{\labelenumi}{\upshape(\theenumi)}

\newtheorem{Theorem}{Theorem}[section]
\newtheorem{Proposition}{Proposition}[section]
\newtheorem{Lemma}{Lemma}[section]

\newenvironment{Proofc}[1]{\smallskip\par\noindent\textsc{#1}\quad}%
  {\hfill$\Box$\bigskip\par}

\newenvironment{Remc}[1]{\smallskip\par\noindent\textsc{#1}\quad}{\smallskip\par}

\newenvironment{myenumerate}{%
\renewcommand{\theenumi}{{\rm(\arabic{enumi})}}%
\renewcommand{\labelenumi}{\theenumi}%
\begin{list}{\labelenumi}%
        {%
        \usecounter{enumi}%
        }}%
{\end{list}}
{\end{myenumerate}}

\def\a{\alpha}

\def\g{\gamma}
\def\l{\lambda}

\def\s{\sigma}

\def\th{\theta}

\def\e{\varepsilon}

\newcommand{\hx}{{\hat x}}

\newcommand{\tx}{{\tilde x}}

\newcommand{\Ho}{{\overline H}}

\newcommand{\xo}{{\overline x}}

\newcommand{\po}{{\overline p}}
\newcommand{\Xo}{{\overline X}}
\newcommand{\Co}{{\overline C}}

\newcommand{\cA}{{\cal A}}

\newcommand{\cF}{{\cal F}}
\newcommand{\cP}{{\cal P}}

\newcommand{\p}{\partial}

\newcommand{\tr}{\operatorname{\text{tr}}}
\newcommand{\N}{{\mathbb N}}
 \newcommand{\Z}{{\mathbb Z}}
\newcommand{\R}{{\mathbb R}}
\newcommand{\Rn}{{\mathbb R^n}}

\newcommand{\Sym}{{\mathbb S}}
\newcommand{\E}{{\mathbb E}}

\newcommand{\esssup}{\mathop{\rm ess{\,}sup}}

\begin{document}
\title{Continuous dependence estimates for the ergodic problem of Bellman equation with an application to the rate of convergence for the homogenization problem\footnotemark[1] }
\author{Claudio Marchi\footnotemark[2]}
\date{version: \today}\maketitle

\footnotetext[1]{Work partially supported by
the INDAM-GNAMPA project ``Fenomeni di propagazione su grafi ed in mezzi eterogenei''.}
\footnotetext[2]{Dipartimento di Matematica, Universit\`a di Padova, via Trieste 63, 35121 Padova, Italy ({\tt marchi@math.unipd.it}).}
 \begin{abstract}
This paper is devoted to establish continuous dependence estimates for the ergodic problem for Bellman operators (namely, estimates of $(v_1-v_2)$ where~$v_1$ and~$v_2$ solve two equations with different coefficients).
We shall obtain an estimate of $\|v_1-v_2\|_\infty$  with an explicit dependence on the $L^\infty$-distance between the coefficients and an explicit characterization of the constants and also, under some regularity conditions, an estimate of $\|v_1-v_2\|_{C^2(\R^n)}$.

Afterwards, the former result will be crucial in the estimate of the rate of convergence for the homogenization of Bellman equations.
In some regular cases, we shall obtain the same rate of convergence established in the monographs~\cite{BLP,JKO} for regular linear problems.

 \end{abstract}

\begin{description}
\item [{\bf MSC 2000}:] 35B27, 35B30, 35J60, 49L25.
 \item [{\bf Keywords}:] Continuous dependence estimates, Hamilton-Jacobi-Bellman equations, viscosity solutions, ergodic problems, homogenization, rate of convergence.
\end{description}


\section{Introduction}\label{intro}

In this paper we study the continuous dependence for the ergodic problem of (Hamilton-Jacobi-)Bellman operators. This property will be crucial (among others) in the estimate of the rate of convergence for the homogenization problem of such operators.

\vskip2mm

In the {\it ergodic problem}, we seek a pair $(v,U)$, with $v\in C(\Rn)$ and $U\in \R$ (that is, $v$ is a real-valued  function while $U$ is a constant) which, in viscosity sense, satisfies
\begin{equation}\label{E1}
H(x,Dv,D^2v)=U\qquad\textrm{in }\R^n
\end{equation}
where the Bellman operator
\begin{equation}\label{max}
H(x,p,X)=\max\limits_{\a\in A}\left\{
-\tr\left(a(x,\a)X\right)-f(x,\a)\cdot p -\ell(x,\a)\right\}
\end{equation}
is uniformly elliptic and periodic in~$x$
(we recall that these operators naturally arise in optimal control problems; see below).
This kind of equations have been widely studied in literature, in particular in connection with homogenization or singular perturbation problems (see~\cite{AB8, AL, CSW, Ev2} 
and references therein), with long-time behaviour of solution to parabolic equations (see~\cite{AL,BS01}) 
with dynamical systems in a torus (see~\cite{AA, CFS}) and with the principal eigenvalue and maximum principle (see~\cite{Ev2, AB8}). 
In order to solve problem~\eqref{E1}, it is expedient to introduce  either the {\it stationary approximated} problems
\begin{equation}\label{approx}
\l v^\l+H(x,Dv^\l,D^2v^\l)=0\qquad \textrm{in }\R^n
\end{equation}
for any $\l\in(0,1)$, or the {\it evolutive approximated} problem
\begin{equation}\label{parab_approx}
\p_t V+H(x,DV,D^2V)=0\qquad \textrm{in } (0,+\infty)\times \R^n,\qquad
V(0,x)=0\qquad \textrm{on } \R^n.
\end{equation}
It is well known (see \cite{AL, AB8}) that: (1) there exists exactly one value $U\in\R$ (called {\it ergodic constant}) such that equation~\eqref{E1} admits a periodic solution (which is unique up to a constant); (2) there exist exactly one periodic solution to~\eqref{approx} and one $x$-periodic solution to~\eqref{parab_approx}; (3) as $\l\to 0^+$ (respectively $t\to+\infty$), $-\l v^\l$ and $v^\l-v^\l(0)$ (resp., $V(t,x)/t$ and $V(t,x)-Ut$) uniformly converge to~$U$ and respectively to a solution~$v$ of~\eqref{E1}.

It is worth to recall (see~\cite{FlS}) that problems~\eqref{approx} and~\eqref{parab_approx} naturally arise in stochastic optimal control problem in a periodic medium: consider the control system for $s>0$
\[dx_s= f(x_s,\a_s) ds +\sqrt 2 \s(x_s,\a_s) dW_s,\qquad x_0=x\]
where $(\Omega,\cF,\cP)$ is a probability space endowed with a right continuous filtration~$(\cF_t)_{0\leq t<+\infty}$ and a $p$-dimensional Brownian motion~$W_t$. The control $\a$ is chosen in the set~$\cA$ of progressively measurable processes with value in the compact set~$A$ for minimizing the {\it cost}
\[
P(x,\a):=\E_x\int_0^{+\infty} \ell(x_s,\a_s)e^{-\l s}ds,\qquad
(\textrm{resp., }P(t,x,\a):=\E_x\int_0^t \ell(x_s,\a_s)ds )
\]
where $\E_x$ denotes the expectation.
Then, the value function~$v^\l(x):=\inf_{\a\in\cA}P(x,\a)$ (resp., $V(t,x):=\inf_{\a\in\cA}P(t,x,\a)$) solves~\eqref{approx} (resp., ~\eqref{parab_approx}) with $a=\s \s^T$.
We deduce that the ergodic constant and a solution to~\eqref{E1} can be written as (for some $c_0\in\R$)
\begin{eqnarray*}
&U=\lim\limits_{\l\to0^+} \l \inf\limits_{\a\in\cA}\E_x\int_0^{+\infty} \ell(x_s,\a_s)e^{-\l s}ds
=\lim\limits_{t\to+\infty}\frac 1t \inf\limits_{\a\in\cA}\E_x\int_0^t \ell(x_s,\a_s)ds&\\
&v=\lim\limits_{\l\to0^+} \l\left[\inf\limits_{\a\in\cA}P(x,\a)-\inf\limits_{\a\in\cA}P(0,\a)\right]
= \lim\limits_{t\to+\infty}\left[\inf\limits_{\a\in\cA} P(t,x,\a)-Ut\right]+c_0.&
\end{eqnarray*}
%

In this paper, we tackle some continuous dependence estimates for the solution~$v$ to~\eqref{E1} with $v(0)=0$; in fact, we shall address the estimate of $\|v_1-v_2\|_\infty$ and of $\|v_1-v_2\|_{C^2(\R^n)}$, where $v_1$ and $v_2$ are solutions to two equations~\eqref{E1} with different coefficients. The  additional requirement $v(0)=0$ is needed to avoid the non-uniqueness of $v$.
We recall that the continuous dependence estimate for the ergodic constant has been established in~\cite{AB1} (see \cite{M12} and \cite{AB8} respectively for a generalization to Bellman-Isaacs operators and to some degenerate elliptic cases). 
Let us recall that such estimates have been widely studied because they play a crucial role in many contexts as error estimates for approximation schemes (see \cite{BJ, DK} and references therein), regularity results (for instance, see \cite{BLM, JK1}) and rate of convergence for vanishing viscosity methods (see \cite{JG, JK1} and references therein).
Unfortunately, the previous literature does not apply to~\eqref{E1}; indeed, as far as we know, all the papers consider either elliptic equations with a strictly positive coefficient of the $0$-th order term or parabolic equations where the constant entering in the estimate diverges as $t\to+\infty$.
In other words, these results apply to~\eqref{approx} or to~\eqref{parab_approx} but the estimates would diverge as $\l\to0^+$ or $t\to+\infty$.
Finally, for a diffusion matrix left unchanged, the continuous dependence estimate in $\mathbb L^\infty$-distance is ascertained in~\cite{CCM11}.

\vskip2mm

In the second part of the paper, we shall study the rate of convergence for the homogenization of  Bellman equations.
Consider the problem
\begin{equation}\label{homog_intro}
u^\e+H\left(x,\frac x\e,Du^\e,D^2u^\e\right)=0\qquad \textrm{in }\R^n,
\end{equation}
where $H$ has the form~\eqref{max} and its coefficients depend on $(x,x/\e,\a)$ and are $x/\e$-periodic.
It is well known  (see \cite{AB1,Ev2}), that as $\e\to0^+$, the solution~$u^\e$ converges locally uniformly to the solution of the {\it effective} problem
\begin{equation*}
u+\Ho\left(x,Du,D^2u\right)=0\qquad \textrm{in }\R^n
\end{equation*}
where the value $\Ho(x,p,X)$ is characterized as the ergodic constant for a suitable Hamiltonian (see Section~\ref{sect:rate} for details). This feature is formally motivated by the the expansion (which goes back to the homogenization theory for linear equations; see \cite{BLP,JKO})
\begin{equation}\label{expansion}
u^\e(x)=u(x)+\e v(x,x/\e)+\e^2 w(x,x/\e)+\dots.
\end{equation}

In the framework of fully nonlinear equations, the study of the rate of convergence for the homogenization problem (namely, an estimate of $\|u-u^\e\|_\infty$) started with the paper~\cite{CI01} for first order problems (for the linear case: see~\cite{BLP,JKO} for the periodic environment and~\cite{Y86,BP99} for the random one).
Camilli and the author~\cite{CM09} tackled this issue for the homogenization of Bellman problem~\eqref{homog_intro}.
Taking advantage of the $C^{2,\th}$-regularity of~$u$, they proved that $u^\e$ converges to $u$ with a rate of order $2\th/(2+\th)$, namely $\|u-u^\e\|_\infty\leq C \e^{2\th/(2+\th)}$.
It is worth to recall that  Caffarelli and Souganidis~\cite{CS10} obtained a rate of order $|\ln \e|^{-1/2}$ for the homogenization of possibly nonconvex Hamiltonians in random media (afterwards, these results have been extended to related investigations: e.g., see~\cite{CM11} for multiscale homogenization, \cite{CCM11} for the simultaneous effect of homogenization and vanishing viscosity, \cite{ACSppt} for stochastic Hamilton-Jacobi equations and~\cite{AP11} for $u/\e$-periodic Hamiltonians).
%
%
%
%

In this paper, we shall improve the result of~\cite{CM09} obtaining a rate of order~$\th$. Let us stress that, for $\th=1$, our estimate has the same order as the ``natural'' one given by the formal expansion~\eqref{expansion} and also as the regular case for linear problem as stated in the monographs by Bensoussan, J.L. Lions and Papanicolaou~\cite[pag. 76]{BLP} and by Jikov, Kozlov and Oleinik~\cite[pag. 30]{JKO}.

\vskip2mm

In conclusion, this paper is devoted to two main purposes.
The former is to establish continuous dependence estimates in the $\mathbb L^\infty$-distance (and in the $C^2$-distance) for the solution~$v$ to~\eqref{E1}-\eqref{max} with an explicit dependence on the $\mathbb L^\infty$-distance between the coefficients and an explicit characterization of the constants.
The latter is to establish an estimate of the rate of convergence for the homogenization problem~\eqref{homog_intro} that generalizes to Bellman operators the estimate given in the monographs~\cite{BLP,JKO} for linear ones.

\vskip2mm

This paper is organized as follows: in the  rest of this Section we set some notations.
Section~\ref{sect:cde} is devoted to the continuous dependence estimates for the ergodic problem.
Section~\ref{sect:rate} concerns the rate of convergence for the homogenization problem.

\vskip5mm

\textbf{Notations: } We define $\mathbb M^{n,p}$ and $\Sym^n$ respectively as the set of
$n\times p$ real matrices and the set of $n\times n$ real symmetric matrices. The latter is endowed with the {\it Euclidean} norm and with the usual order, namely: for $X=(X_{ij})_{i,j=1,\dots,n}\in \Sym^n$, $|X|:=\left(\sum_{i,j=1}^n X^2_{ij}\right)^{1/2}\equiv \tr({}^T\!AA)$ and, for $X,Y\in\Sym^n$, we shall write $X\geq Y$, if $X-Y$ is a semi-definite positive matrix.
We denote $I$ the identity matrix in $\Sym^n$.


For every real-valued function $h$, we set $\|h\|_\infty:=\esssup|h(y)|$.
For $\th\in(0,1]$, we use the $\th$-H\"older seminorm: $[h]_{\th,A}:=\sup\{\frac{|h(y)-h(x)|}{|y-x|^\th}\mid y,x\in A,\, y\neq x\}$. $C^\th(\R^n)$ denotes the space of functions $h$ such that: $\|h\|_\infty+[h]_\th<+\infty$ while, for $m \in\N$, $C^{m,\th}(\R^n)$ stands for the space of functions whose derivatives of order $m$ belong to $C^\th(\R^n)$.

For a metric space~$A$ with $x\in A$ and $r>0$, $B_A(x,r)$ stands for the open ball of radius~$r$ centered in $x$, i.e. $B_A(x,r):=\{y\in A\mid \textrm{dist}(x,y)<r\}$.
For every $a,b\in\R$, we set $a\wedge b:=\min\{a,b\}$ and $a\lor b:=\max\{a,b\}$.


%
%
%
%
\section{Continuous dependence estimates for the ergodic problem}\label{sect:cde}

This Section is devoted to estimate the distance, both in the $\mathbb L^\infty$-norm and in the $C^2$-norm, between the functions~$v_1$ and $v_2$ where (for $i=1,2$) $v_i$ is the solution to the following ergodic problem
\begin{equation}\label{Ei}\tag{Ei}
H_i\left(x, Dv_i,D^2v_i\right)=U_i\quad \textrm{in }\R^n,\quad v_i(0)=0
\end{equation}
with
\begin{equation}\label{Hik}
H_i(x,p,X):=\max\limits_{\a\in A}\left\{-\tr\left(a_i(x,\a)X\right)-f_i(x,\a)\cdot p -\ell_i(x,\a)\right\}.
\end{equation}

We shall assume
\begin{itemize}
\item[(A1)] $A$ is a compact metric space;
\item[(A2)] for $\s_i(x,\a)\in\mathbb M^{n,p}$, $a_i=\s_i \s_i^T$ is {\it uniformly elliptic}: for some $\nu>0$, there holds $a(x,\a)\geq \nu I$,  $\forall(x,\a)\in \R^n\times A$;
\item[(A3)] the functions $a_i$, $f_i$ and $\ell_i$ are $\Z^n$-periodic in $x$; moreover, they are bounded and Lipschitz continuous in~$x$ uniformly in~$\a$, namely, for $\phi=a_i,f_i,\ell_i$, there holds
\[\|\phi\|_\infty\leq K_\phi,\qquad |\phi(x,\a)-\phi(y,\a)|\leq L_\phi|x-y|\quad \forall x,y\in\R^n,\, \a\in A.\]
\end{itemize}

\vskip2mm

In order to study problem~\eqref{Ei}, it is expedient to introduce the approximate equations for $\l\in(0,1)$
\begin{equation}\label{Ai}\tag{Ai}
\l v^\l_i+H_i\left(x,Dv^\l_i,D^2v^\l_i\right)
=0\qquad \textrm{in }\R^n.
\end{equation}
In the following Proposition, we present useful properties of solutions of problems~\eqref{Ei} and~\eqref{Ai}. Though most of them are already known in the literature, we need to be careful for the details of dependence for constants which appear in some estimate for the later argument. Therefore, we give the proof at the end of this Section for reader's convenience.

\begin{Proposition}\label{l1p2}
The following properties hold
\begin{itemize}
\item[(i)] There exists a unique periodic solution to~\eqref{Ai} with~$\|\l v_i^\l\|_\infty\leq K_{\ell_i}$.
\item[(ii)] For some constant~$\th\in(0,1]$ depending only on $n$ and $\nu$, there holds
\[\|v_i^\l-v_i^\l(0)\|_{C^{2,\th}(\R^n)}\leq C_1(1+K_{\ell_i}+L_{\ell_i})\qquad \forall \l\in(0,1)\]
where $C_1$ is a constant independent of $\l$, $K_{\ell_i}$ and $L_{\ell_i}$.
\item[(iii)] As $\l\to 0^+$, the sequence~$\{\l v_i^\l\}_\l$ uniformly converges to a value~$-U_i$; moreover, $U_i$ is the unique constant such that the problem~\eqref{Ei} admits a periodic solution.
\item[(iv)] As $\l\to 0^+$, the sequence~$\{v_i^\l-v_i^\l(0)\}_\l$ converges in the $C^2(\R^n)$-norm to the unique solution~$v_i$ to~\eqref{Ei}.
Furthermore, there holds
\[\|v_i\|_{C^{2,\th}(\R^n)}\leq C_1(1+K_{\ell_i}+L_{\ell_i})\]
where $\th$ and $C_1$ are the constants introduced in point $(ii)$.
\end{itemize}
\end{Proposition}

Let us now establish our main result on the $\mathbb L^\infty$-continuous dependence estimate for problem~\eqref{Ei}.

\begin{Theorem}\label{thm:cde}
Let $v_i$ be the unique solution to problem~\eqref{Ei} ($i=1,2$). Then, there holds
\begin{equation*}
\|v_1-v_2\|_\infty \leq C C_\ell\left(\max_{x,\a}|a_1-a_2|+\max_{x,\a}|f_1-f_2|\right)+ C\max_{x,\a}|\ell_1-\ell_2|
\end{equation*}
where $C_\ell:=1+(K_{\ell_1}+L_{\ell_1})\wedge(K_{\ell_2}+L_{\ell_2})$ while $C$ is a constant independent of $K_{\ell_i}$ and $L_{\ell_i}$.
\end{Theorem}
\begin{Proofc}{Proof of Theorem~\ref{thm:cde}}
For $i=1,2$, let $v^\l_i$ be the solution to the approximated problem~\eqref{Ai}; set  $w^\l_i:=v^\l_i-v^\l_i(0)$. By Proposition~\ref{l1p2}-$(iv)$, it suffices to establish
\begin{equation}\label{s3p3}
\|w^\l_1-w^\l_2\|_\infty \leq C C_\ell\left(\max_{x,\a}|a_1-a_2|+\max_{x,\a}|f_1-f_2|\right)+ C\max_{x,\a}|\ell_1-\ell_2|\qquad\forall \l\in(0,1)
\end{equation}
where $C$ is a constant depending only on $K_{a_i}$, $L_{a_i}$, $K_{f_i}$, $L_{f_i}$, $n$ and $\nu$ (i.e., it is independent of~$K_{\ell_i}$, $L_{\ell_i}$ and $\l$).

Let us claim that, for every $\l\in(0,1)$, there holds
\begin{equation}\label{s1p4}
\l \|v^\l_1-v^\l_2\|_\infty \leq C_1 C_\ell\left(\max_{x,\a}|a_1-a_2|+\max_{x,\a}|f_1-f_2|\right)+ \max_{x,\a}|\ell_1-\ell_2|,
\end{equation}
where $C_1$ and $C_\ell$ are the constants introduced respectively in Proposition~\ref{l1p2}-$(ii)$ and in the statement.
Indeed, without any loss of generality one can assume $K_{\ell_1}+L_{\ell_1}\geq K_{\ell_2}+L_{\ell_2}$ and easily check that the functions
\[v^\pm:=v^\l_2\pm \l^{-1}C_1 C_\ell\left(\max_{x,\a}|a_1-a_2|+\max_{x,\a}|f_1-f_2|\right)+ \l^{-1}\max_{x,\a}|\ell_1-\ell_2|\]
are respectively a super- and a subsolution to problem~\eqref{Ai} with $i=1$.
The comparison principle guarantees $v^-\leq v^\l_1\leq v^+$ which is equivalent to~\eqref{s1p4}.

In order to prove inequality~\eqref{s3p3}, we shall proceed by contradiction assuming that, for $k\in \N$ and $i=1,2$, there exist $\l_k$, $a_{ik}$, $f_{ik}$ and $\ell_{ik}$ such that
\begin{itemize}
\item[-)] $\l_k\to 0$ as $k\to+\infty$;
\item[-)] $a_{ik}$ and $f_{ik}$ satisfy $(A2)$-$(A3)$ with the same constants $\nu$, $K_a$, $L_a$, $K_f$ and $L_f$ while $\ell_{ik}$ satisfies $(A3)$ with some constants~$K_{\ell_{ik}}$ and~$L_{\ell_{ik}}$;
\item[-)] denote $v^{\l_k}_i$ the solution to~\eqref{Ai} with $\l$, $a_i$, $f_i$ and $\ell_i$ replaced respectively by $\l_k$, $a_{ik}$, $f_{ik}$ and $\ell_{ik}$;
the functions $w^{\l_k}_i:=v^{\l_k}_i-v^{\l_k}_i(0)$ verify
\begin{equation*}
c_k:=\|w^{\l_k}_1-w^{\l_k}_2\|_\infty\geq k [C_{k}\left(\max_{x,\a}|a_{1k}-a_{2k}|+\max_{x,\a}|f_{1k}-f_{2k}|\right)+ \max_{x,\a}|\ell_{1k}-\ell_{2k}|]
\end{equation*}
with $C_k:=1+(K_{\ell_{1k}}+L_{\ell_{1k}})\wedge(K_{\ell_{2k}}+L_{\ell_{2k}})$. 
\end{itemize}
Without any loss of generality, we assume $K_{\ell_{1k}}+L_{\ell_{1k}} \geq K_{\ell_{2k}}+L_{\ell_{2k}}$ for every $k\in \N$; hence, we have $C_k=1+ K_{\ell_{2k}}+L_{\ell_{2k}}$.
Taking advantage of the inequality $\max_A f-\max_A g\leq \max_A(f-g)$ for every $f,g\in C(A)$, we observe that the function $\tilde w^k:=w^{\l_k}_1-w^{\l_k}_2$ solves
\begin{multline*}
\l_k \tilde w^k +\l_k\left(v^{\l_k}_1(0)-v^{\l_k}_2(0)\right) 
+\max_{\a\in A}\left\{-\tr(a_{1k} D^2 \tilde w^k)-f_{1k}\cdot D\tilde w^k\right\}\\
+\max_{\a\in A}\left\{\tr[(a_{2k}-a_{1k}) D^2 w^{\l_k}_2]+(f_{2k}-f_{1k})\cdot D w^{\l_k}_2+\ell_{2k}-\ell_{1k}\right\}\geq 0.
\end{multline*}
In particular, the function~$w_k:=\tilde w ^k/c_k\equiv \tilde w ^k/\|\tilde w ^k\|_\infty$ fulfills
\begin{multline*}
H^k(x,D  w_k,D^2 w_k)+\l_k  w_k +c_k^{-1}\l_k \left(v^{\l_k}_1(0)-v^{\l_k}_2(0)\right) \\ +c_k^{-1}\max_{\a\in A}\left\{\tr[(a_{2k}-a_{1k}) D^2 w^{\l_k}_2]+(f_{2k}-f_{1k})\cdot D w^{\l_k}_2+\ell_{2k}-\ell_{1k}\right\}\geq 0,
\end{multline*}
with 
\[H^k(x,p,X):=\max_{\a\in A}\left\{-\tr(a_{1k}(x,\a) X)-f_{1k}(x,\a)\cdot p\right\}.\]
Since $\|w_k\|_\infty=1$, by~\eqref{s1p4} and the definition of~$c_k$, we infer
\[
\l_k  w_k +c_k^{-1}\l_k \left(v^{\l_k}_1(0)-v^{\l_k}_2(0)\right)=o(1)\qquad \textrm{as }k\to+\infty.
\]
Moreover, owing to Proposition~\ref{l1p2}-$(ii)$ and to our definition of $c_k$ we have
\begin{align*}
&c_k^{-1}\left|\max_{\a\in A}\left\{\tr[(a_{2k}-a_{1k}) D^2 w^{\l_k}_2]+(f_{2k}-f_{1k})\cdot D w^{\l_k}_2+\ell_{2k}-\ell_{1k}\right\}\right|\\
&\qquad\leq c_k^{-1}
\left(\max_{x,\a}|a_{1k}-a_{2k}|\|D^2 w^{\l_k}_2\|_\infty+\max_{x,\a}|f_{1k}-f_{2k}|\|D w^{\l_k}_2\|_\infty+ \max_{x,\a}|\ell_{1k}-\ell_{2k}|\right)\\
&\qquad\leq c_k^{-1}
\left[C_1(1+K_{\ell_{2k}}+L_{\ell_{2k}})\left(\max_{x,\a}|a_{1k}-a_{2k}|+\max_{x,\a}|f_{1k}-f_{2k}|\right)+ \max_{x,\a}|\ell_{1k}-\ell_{2k}|\right]\\
&\qquad\leq o(1)\qquad \textrm{as }k\to+\infty.
\end{align*}
Taking into account the last three relations, we obtain
\begin{equation}\label{s1p5}
H^k(x,D  w_k,D^2 w_k)\geq o(1)\qquad \text{as }k\to+\infty.
\end{equation}

On the other hand, taking advantage of the inequality $\max_A f-\max_A g\geq \min_A(f-g)$,  we obtain that the function~$w_k$ also fulfills
\begin{multline*}
h^k(x,D w_k,D^2 w_k)+\l_k  w_k +c_k^{-1}\l_k \left(v^{\l_k}_1(0)-v^{\l_k}_2(0)\right) \\ +c_k^{-1}\min_{\a\in A}\left\{\tr[(a_{2k}-a_{1k}) D^2 w^{\l_k}_2]+(f_{2k}-f_{1k})\cdot D w^{\l_k}_2+\ell_{2k}-\ell_{1k}\right\}\leq 0,
\end{multline*}
with $h^k(x,p,X):=\min_{\a\in A}\left\{-\tr(a_{1k}(x,\a) X)-f_{1k}(x,\a)\cdot p\right\}$. Arguing as before, we infer
\begin{equation}\label{s1p5bis}
h^k(x,D  w_k,D^2 w_k)\leq o(1)\qquad \text{as }k\to+\infty.
\end{equation}
Owing to relations~\eqref{s1p5} and~\eqref{s1p5bis}, standard regularity theory for elliptic equations (see \cite{t1,t2}) ensures that the family $\{w_k\}_k$ is uniformly H\"older continuous. Indeed, let us briefly show how to apply the arguments of \cite[Theorem 2.1]{t2} for the H\"older continuity of~$w_k$. Being a subsolution to~\eqref{s1p5bis}, the function~$w_k$ verifies a weak Harnack inequality. On the other hand, since~$w_k$ is a supersolution to~\eqref{s1p5}, also the function~$-w_k$ verifies another weak Harnack inequality.
We observe that in these two weak Harnack inequalities the occurring constants depend only on $n$, $\nu$, $K_a$, $L_a$, $K_f$ and $L_f$ (recall that $\|w_k\|_\infty =1$).
Finally, these two inequalities yield that the function~$w_k$ is H\"older continuous with an H\"older exponent depending only on $n$, $\nu$, $K_a$ and an H\"older constant depending only on $n$, $\nu$, $K_a$, $L_a$, $K_f$ and $L_f$. Therefore, the proof that the family~$\{w_k\}_k$ is uniformly H\"older continuous is accomplished.

We note that~$H^k$ fulfills the following properties for every $(x,p,X)$ and $k\in\N$
\begin{align*}
|H^k(x,p,X)|&\leq K_a|X|+K_f|p|\\
|H^k(x_1,p_1,X_1)-H^k(x_2,p_2,X_2)|&\leq
K_a|X_1-X_2|+K_f|p_1-p_2|\\
&\qquad+[L_a(|X_1|\wedge |X_2|)+L_f(|p_1|\wedge |p_2|)]|x_1-x_2|.
\end{align*}
By the Ascoli Theorem, as $k\to+\infty$ (passing to a subsequence, if necessary), we can assume that $H^k$ converges to an uniformly elliptic operator~$\tilde H$ locally uniformly in $\R^n\times\R^n\times\Sym^n$ and that the function~$w_k$ converges to some continuous periodic function~$w$ uniformly in $\R^n$. Finally, passing to the limit in~\eqref{s1p5}, the stability result for viscosity solution entails that~$w$ is a solution to
\[\tilde H(x,D w,D^2 w)\geq 0.\]
We note that $w$ also verifies: $\|w\|_\infty=1$, and $w(0)=0$; thus, it attains its global minimum. This fact contradicts the strong maximum principle (see: \cite{t1, BL99}); the proof of~\eqref{s3p3} is accomplished.
\end{Proofc}

Let us now establish our main result on the $C^2$-continuous dependence estimate for problem~\eqref{Ei}. To this end it is expedient to set
\begin{equation}\label{def:calK}
 R:=C_1(1+K_\ell+L_\ell)\qquad\textrm{and}\qquad
{\cal K}:=\R^n\times B_{\R^n}(0,R)\times B_{\Sym^n}(0,R)
\end{equation}
where $C_1$ is the constant introduced in Proposition~\ref{l1p2}.

\begin{Theorem}\label{thm:cdecont}
Assume that the coefficients of operators $H_1$ and $H_2$ satisfy hypotheses $(A2)$-$(A3)$ with the same constants $\nu$,  $K_{\phi}$ and $L_{\phi}$, for $\phi=a,f,\ell$.
Assume that, for some constants $K_H$ and $\th'\in(0,1]$, there holds
$\|H_i\|_{C^{1,\th'}({\cal K})}\leq K_H$ for $i=1,2$.
Let $v_i$ be the unique solution to problem~\eqref{Ei} ($i=1,2$). Then, we have
\begin{equation*}
\|v_1-v_2\|_{C^{2}(\R^n)} \leq C \left(\max_{x,\a}|a_1-a_2|+\max_{x,\a}|f_1-f_2|+\max_{x,\a}|\ell_1-\ell_2|\right)+[H_1-H_2]_{1,{\cal K}}
\end{equation*}
where $C$ is a constant depending only on $n$, $K_H$, $\th'$, $\nu$, $K_{\phi}$ and $L_{\phi}$, for $\phi=a,f,\ell$.
\end{Theorem}
\begin{Proofc}{Proof of Theorem~\ref{thm:cdecont}}
We shall adapt the arguments of the proof of Theorem~\ref{thm:cde}; hence, we shall only emphasize the main differences.
As before, for $i=1,2$, we set $w^\l_i:=v^\l_i-v^\l_i(0)$  (recall that $v^\l_i$ is the solution to~\eqref{Ai}) and we observe that estimate~\eqref{s1p4} still holds. By Proposition~\ref{l1p2}-$(iv)$, it suffices to establish that there exists a constant~$C$ such that, for every $\l\in(0,1)$, there holds
\begin{equation}\label{2.1bis}
\|w^\l_1-w^\l_2\|_{C^{2}(\R^n)} \leq C\left(\max_{x,\a}|a_1-a_2|+\max_{x,\a}|f_1-f_2|+ \max_{x,\a}|\ell_1-\ell_2|\right)+[H_1-H_2]_{1,{\cal K}}.
\end{equation}
In order to prove this inequality, we shall proceed by contradiction assuming that, for each $k\in \N$ and $i=1,2$, there exist $\l_k$, $a_{ik}$, $f_{ik}$ and $\ell_{ik}$ such that
\begin{itemize}
\item[-)] $\l_k\to 0$ as $k\to+\infty$;
\item[-)] $a_{ik}$,  $f_{ik}$ and $\ell_{ik}$ satisfy $(A2)$-$(A3)$ with the same constants $\nu$, $K_\phi$ and $L_\phi$ ($\phi=a,f,\ell$);
\item[-)] let $H_{ik}$ be the operator obtained replacing $a_i$, $f_i$ and $\ell_i$ respectively with $a_{ik}$, $f_{ik}$ and $\ell_{ik}$ in~\eqref{Hik}: $\|H_{ik}\|_{C^{1,\th'}({\cal K})}\leq K_H$;
\item[-)] let $v^{\l_k}_i$ solve~\eqref{Ai} with $\l$, and $H_i$ replaced respectively by $\l_k$, and $H_{ik}$: for $w^{\l_k}_i:=v^{\l_k}_i-v^{\l_k}_i(0)$, the value $c_k:=\|w^{\l_k}_1-w^{\l_k}_2\|_{C^{2}(\R^n)}$ verifies
\begin{equation}\label{2.2ter}
c_k\geq k \left(\max_{x,\a}|a_{1k}-a_{2k}|+\max_{x,\a}|f_{1k}-f_{2k}|+ \max_{x,\a}|\ell_{1k}-\ell_{2k}|\right)+[H_{1k}-H_{2k}]_{1,{\cal K}}.
\end{equation}
\end{itemize}
For the sake of simplicity, let us introduce the notation $H_{ik}[\psi]:=H_{ik}(x,D\psi(x),D^2\psi(x))$ for every $\psi\in C^2(\R^n)$.

The function $w_k:=(w^{\l_k}_1-w^{\l_k}_2)/c_k\equiv (w^{\l_k}_1-w^{\l_k}_2)/\|w^{\l_k}_1-w^{\l_k}_2\|_{C^{2}(\R^n)}$ fulfills
\begin{equation}\label{B1}
R_k(x)+c_k^{-1}\left(H_{1k}[w^{\l_k}_1]-H_{1k}[w^{\l_k}_2]\right)=0
\end{equation}
with $R_k:=\l_k w_k +c_k^{-1}\l_k(v^{\l_k}_1(0)-v^{\l_k}_2(0)) +c_k^{-1}(H_{1k}[w^{\l_k}_2]-H_{2k}[w^{\l_k}_2])$.
Our regularity assumption on~$H_{ik}$ yields
\begin{equation}\label{D3}
c_k^{-1}\left(H_{1k}[w^{\l_k}_1]-H_{1k}[w^{\l_k}_2]\right)=
f^k(x)\cdot Dw_k+H^k(x)\cdot D^2w_k
\end{equation}
where
\begin{eqnarray*}
f^k(x)&:=&\int_0^1D_p H_{ik}\left(x,t Dw^{\l_k}_1+(1-t)Dw^{\l_k}_2,t D^2w^{\l_k}_1+(1-t)D^2w^{\l_k}_2\right)\,dt,\\
H^k(x)&:=&\int_0^1D_X H_{ik}\left(x,t Dw^{\l_k}_1+(1-t)Dw^{\l_k}_2,t D^2w^{\l_k}_1+(1-t)D^2w^{\l_k}_2\right)\,dt.
\end{eqnarray*}
Substituting relation~\eqref{D3} in~\eqref{B1}, we get
\begin{equation}\label{E11}
R_k(x)+f^k(x)\cdot Dw_k+H^k(x)\cdot D^2w_k=0.
\end{equation}
Moreover, Proposition~\ref{l1p2}-($ii$) entails that $w^{\l_k}_i$ are equibounded in $C^{2,\th}(\R^n)$; in particular, we deduce that $f^k$ and $H^k$ are uniformly bounded and uniformly H\"older continuous (with exponent $\th\th'$).
Let us now claim that
\begin{equation}\label{C1}
[R_k]_{\th,\R^n}\leq 2R \qquad \forall k\in\N\qquad\textrm{and}\qquad\|R_k\|_\infty=o(1)\quad \textrm{as }k\to+\infty,
\end{equation}
where $\th$ and~$R$ are the constants introduced in Proposition~\ref{l1p2}-($ii$) and respectively in~\eqref{def:calK}.
Actually, from Proposition~\ref{l1p2}-($ii$) we recall that $\|w^{\l_k}_2\|_{C^{2,\th}(\R^n)}\leq R$; hence, we deduce
\[
c_k^{-1}[H_{1k}[w^{\l_k}_2]-H_{2k}[w^{\l_k}_2]]_{\th,\R^n}\leq
c_k^{-1}[H_{1k}-H_{2k}]_{1,{\cal K}}\|w^{\l_k}_2\|_{C^{2,\th}(\R^n)}\leq R. 
\]
Taking into account relation $\|w_k\|_{C^{2}(\R^n)}=1$, we infer the first part of~\eqref{C1}. Moreover, by the same arguments as those in the proof of Theorem~\ref{thm:cde}, we achieve the proof of the second part of~\eqref{C1}.

Invoking standard regularity theory for linear elliptic equations (see~\cite{GT,LU}), we deduce that, for some $\th''\in(0,1]$, the functions~$w_k$ are uniformly bounded in $C^{2,\th''}(\R^n)$.
By the Ascoli theorem (possibly passing to a subsequence), $w_k$ uniformly converges to some periodic function~$w\in C^{2,\th''}(\R^n)$ along with all its derivatives up to order~$2$ while $f^k$ and $H^k$ converge locally uniformly to some functions~$\tilde f$ and~$\tilde H$ respectively. Passing to the limit as $k\to+\infty$ in~\eqref{E11}, by~\eqref{C1}, the stability result ensures that $w$ is a solution
\begin{equation*}
\tilde f(x)\cdot Dw+\tilde H(x)\cdot D^2w=0.
\end{equation*}
The strong maximum principle yields that, being periodic, $w$ must be constant; furthermore, since $w_k(0)=0$, we get: $w\equiv 0$. On the other hand, passing to the limit in $\|w_k\|_{C^{2}(\R^n)}=1$, we deduce $\|w\|_{C^{2}(\R^n)}=1$ which gives the desired contradiction. Whence, the proof of estimate~\eqref{2.1bis} is accomplished.
\end{Proofc}

\begin{Proofc}{Proof of Proposition~\ref{l1p2}} For the detailed proof, we refer the reader to~\cite[Theorem II.2]{AL} and to~\cite[Proposition 12]{AB1} (see also~\cite[Theorem 4.1]{AB8} for similar results for Bellman-Isaacs operators).
Let us just discuss the proof of point~$(ii)$ and, especially, its right hand side which seems to be new in this form.
For the sake of simplicity, we shall drop the subscript {\it ``i''}.
We claim that there exists a constant~$C'_1$, enjoying the same properties of~$C_1$ (i.e., independent of $\l$, $K_\ell$ and $L_\ell$) such that
\begin{equation}\label{1p2}
\|v^\l-v^\l(0)\|_\infty\leq C'_1(1+K_\ell)\qquad \forall \l\in(0,1).
\end{equation}
In order to prove this estimate, we proceed by contradiction assuming that, for $k\in\N$, there exist $\l_k$, $a_k$, $f_k$ and $\ell_k$, fulfilling~$(A1)$-$(A3)$ with the same constants $K_a$, $L_a$, $K_f$ and~$L_f$, such that $\l_k\to 0$ as $k\to +\infty$ and
\begin{equation}\label{2p2}
\eta_k:=\|v^{\l_k}-v^{\l_k}(0)\|_\infty\geq k(1+K_{\ell_k})\qquad \forall k\in\N
\end{equation}
where~$v^{\l_k}$ is the solution to~\eqref{Ai} with $\l$, $a_i$, $f_i$ and $\ell_i$ replaced respectively by $\l_k$, $a_k$, $f_k$ and $\ell_k$.

The function~$v^k:=\eta_k^{-1}(v^{\l_k}-v^{\l_k}(0))$ fulfills
\[\l_k v^k+\eta_k^{-1}\l_k v^{\l_k}(0)+\max_{\a}
\left\{-\tr\left(a_k D^2v^k\right)-f_k\cdot Dv^k -\eta_k^{-1}\ell_k\right\}=0.\]
By $\|v^k\|_\infty=1$, relation~\eqref{2p2} and point~($i$), we have
\[ \|\l_k v^k\|_\infty +|\eta_k^{-1}\l_k v^{\l_k}(0)|+\max_{\a,x}|\eta_k^{-1}\ell_k|=o(1)\qquad \textrm{as }k\to+\infty.\]
The last two relations entail
\[\max_{\a}\left\{-\tr\left(a_k D^2v^k\right)-f_k\cdot Dv^k\right\}=o(1)\qquad \textrm{as }k\to+\infty.\]
The standard regularity theory for uniformly elliptic HJB equations (see~\cite{GT,Sa88}) guarantees that the family~$\{v^k\}_k$ is uniformly H\"older continuous. Passing to a subsequence if necessary, we assume that~$v^k$ uniformly converges to some periodic function~$\tilde v$ with $\tilde v(0)=0$ and $\|\tilde v\|_\infty=1$.
Moreover, for any $\a_0\in A$ fixed, the previous inequality gives 
\[
-\tr\left(a_k(\cdot,\a_0) D^2v^k\right)-f_k(\cdot,\a_0)\cdot Dv^k
\leq o(1) \qquad\textrm{as }k\to+\infty.\]
By $(A3)$ (passing to a subsequence, if necessary), $a_k(\cdot,\a_0)$ and $f_k(\cdot,\a_0)$ are uniformly convergent respectively to some $\tilde a$ and $\tilde f$. 
Letting $k\to +\infty$ in the previous inequality, by point~$(i)$ and relation~\eqref{2p2}, the stability result of viscosity solutions ensures
\[-\tr\left(\tilde a D^2\tilde v\right)-\tilde f \cdot D\tilde v \leq 0.\]
By the strong maximum principle (see~\cite{t1,BL99}), $\tilde v$ cannot attain its maximum; this property contradicts the fact that $\tilde v$ is periodic with $\tilde v(0)=0$ and $\|\tilde v\|_\infty=1$. hence, our claim~\eqref{1p2} is established.

Finally, invoking~\cite[Theorem 1.1]{Sa88}, for some $\th\in (0,1]$ depending on $n$ and $\nu$ and for some constant $C''_1$ enjoying the same properties of $C_1$, we get
\[\|v^\l-v^\l(0)\|_{C^{2,\th}(\R^n)}\leq C''_1\left(C'_1(1+K_\ell)+L_\ell\right)\]
which amounts to our statement.
\end{Proofc}

%
%
%
%
\section{Rate of convergence for the homogenization problem}\label{sect:rate}

We consider the following homogenization problem
\begin{equation}\label{homog}
u^\e+H\left(x,\frac x\e,Du^\e,D^2u^\e\right)=0\qquad \textrm{in }\R^n
\end{equation}
with
\begin{equation}\label{Fmax}
H(x,y,p,X)=\max\limits_{\a\in A}\left\{
-\tr\left(a(x,y,\a)X\right)-f(x,y,\a)\cdot p -\ell(x,y,\a)\right\}.
\end{equation}
Throughout this Section, the following assumptions will hold
\begin{itemize}
\item[(H1)] $A$ is a compact metric space;
\item[(H2)] for $\s(x,y,\a)\in\mathbb M^{n,p}$, $a=\s \s^T$ is uniformly elliptic: for some $\nu>0$, there holds $a(x,y,\a)\geq \nu I$,  $\forall(x,y,\a)\in \R^n\times\R^n\times A$;
\item[(H3)] the functions $a$, $f$ and $\ell$ are $\Z^n$-periodic in~$y$; moreover, they are bounded and Lipschitz continuous in~$(x,y)$ uniformly in~$\a$: for $\phi=a,f,\ell$, there holds $\|\phi\|_\infty\leq K$ and 
\[ |\phi(x_1,y_1,\a)-\phi(x_2,y_2,\a)|\leq L(|x_1-x_2|+|y_1-y_2|)\quad \forall x_1,x_2,y_1,y_2\in\R^n, \a\in A.\]
\end{itemize}
Let us recall that equation~\eqref{homog} arises, for instance, in stochastic optimal control problem in a medium with microscopic periodic heterogeneities. Actually, consider the dynamics
\[dx_s= f(x_s,x_s/\e,\a_s) ds +\sqrt 2 \s(x_s,x_s/\e,\a_s) dW_s,\qquad x_0=x\]
in the probability space~$(\Omega,\cF,\cP)$ with a right continuous filtration~$(\cF_t)_{0\leq t<+\infty}$ and a $p$-dimensional Brownian motion~$W_t$.
It is well known (see \cite{FlS}) that there holds
\[u^\e(x)=\inf_{\a\in\cA}\E_x\int_0^{+\infty} \ell(x_s,x_s/\e,\a_s)e^{-s}ds\]
where~$\E_x$ denotes the expectation while $\cA$ stands for the set of progressively measurable processes with value in~$A$.
The aim of homogenization is to investigate the behaviour of the optimal control problem as the heterogeneities become smaller and smaller. It is well known (see~\cite{Ev2,AB1}; we refer the reader to~\cite{AB8} for more general classes of operators) that $u^\e$ converges locally uniformly to a limit function~$u$ which can be characterized as the solution of a suitable problem (the {\it effective problem}).
This Section is devoted to improve the estimate of $\|u^\e-u\|_\infty$ established in~\cite{CM09}.

As in~\cite{Ev2,AB1}, we introduce the {\it effective} operator~$\Ho$ as follows: for every $(\xo,\po,\Xo)\in\R^n\times\R^n\times\Sym^n$ fixed, the value $\Ho(\xo,\po,\Xo)$ is the ergodic constant for $H(\xo,\cdot, \po, \Xo+\cdot)$, namely it is the unique constant such that there exists a solution to the {\it cell problem}
\begin{equation}\label{s1p10}
H(\xo,y,\po,\Xo+D^2v)=\Ho(\xo,\po,\Xo).
\end{equation}
The effective problem is
\begin{equation}\label{eff}
u+\Ho\left(x,Du,D^2u\right)=0\qquad \textrm{in }\R^n.
\end{equation}
In the next two statement we shall collect some properties of problems~\eqref{homog} and~\eqref{eff} and respectively of their solutions. We shall omit the proofs: the proof of the former can be found in~\cite[Proposition 12]{AB1} and in~\cite[Lemma 3.2]{Ev2} (see also~\cite[Lemma 2.2]{CM09}).
The proof of the latter is an easy consequence of the previous one and of the regularity theory for uniformly elliptic convex operators (see~\cite[Theorem 1.1]{Sa88}).

\begin{Lemma}\label{l2p9}
The effective operator~$\Ho$ is convex in~$X$ and uniformly elliptic with the same ellipticity constant~$\nu$ as in $(H2)$. Moreover, there exists a constant~$\Co$ such that
\begin{multline*}
|\Ho(x_1,p_1,X_1)-\Ho(x_2,p_2,X_2)|\leq \Co(1+|p_1|\wedge|p_2|+|X_1|\wedge|X_2|)|x_1-x_2|\\+\Co(|p_1-p_2|+|X_1-X_2|)
\end{multline*}
for every $(x_i,p_i,X_i)\in\R^n\times\R^n\times \Sym^n$, $i=1,2$.
Hence, the comparison principle holds for problem~\eqref{eff}.
\end{Lemma}

\begin{Lemma}\label{l3}
The problems~\eqref{homog} and~\eqref{eff} admit exactly one bounded viscosity solution~$u^\e$ and respectively~$u$.
Moreover, there exist $N>0$ and $\th\in(0,1]$ such that
\[\|u^\e\|_\infty\leq N,\quad
\|u\|_\infty\leq N,\quad \|Du\|_\infty\leq N,\quad
\|u\|_{C^{2,\th}(B(x,1))}\leq N \quad \forall x\in\R^n.
\]
\end{Lemma}

Let us now establish our main result on the rate of convergence for the homogenization problem~\eqref{homog}.

\begin{Theorem}\label{thm:hom} 
Let $u^\e$ and $u$ be respectively the solution to problems~\eqref{homog} and~\eqref{eff}. Then, there exists a constant $M$ such that
\[\|u^\e-u\|_\infty\leq M\e^\th\qquad \forall\e\in(0,1)
\]
where $\th\in(0,1]$ is the constant introduced in Lemma~\ref{l3}.
\end{Theorem}
\begin{Proofc}{Proof of Theorem~\ref{thm:hom}}
We shall proceed using some techniques introduced in~\cite[Theorem 2.1]{CM09}. 
We shall denote $v(\cdot;\xo,\po,\Xo)$ the unique solution to~\eqref{s1p10} with $v(0,\xo,\po,\Xo)=0$ so as to display its dependence on the fixed parameters. In the following statement, we collect some regularity properties of~$v$.

\begin{Lemma}\label{cor:p11}
There exist $\th\in(0,1]$ and a constant~$C_2$ such that
\begin{align*}
(i)&\qquad \|v(\cdot;x_1,p_1,X_1)\|_{C^{2,\th}(\R^n)}\leq C_2(1+|p_1|+|X_1|),\\
(ii)&\qquad \|v(\cdot;x_1,p_1,X_1)-v(\cdot;x_2,p_2,X_2)\|_\infty
\leq C_2(|p_1-p_2|+|X_1-X_2|)
\\&\qquad \qquad\qquad +C_2[1+(|p_1|+|X_1|)\wedge(|p_2|+|X_2|)]|x_1-x_2|,
\end{align*}
for every $x_1,x_2,p_1,p_2\in\R^n$ and $X_1,X_2\in \Sym^n$.
\end{Lemma}
\begin{Proofc}{Proof of Lemma~\ref{cor:p11}}
We note that equation~\eqref{s1p10} fulfills assumptions~$(A1)$-$(A3)$ with constants $K_a=K_f=K$, $L_a=L_f=L$, $K_\ell=K(1+|\po|+|\Xo|)$ and $L_\ell=L(1+|\po|+|\Xo|)$. Whence, point~$(i)$ is due to Proposition~\ref{l1p2}-$(ii)$ and -$(iv)$.

In order to ascertain point~$(ii)$, it suffices to apply Theorem~\ref{thm:cde} with the variable~$x$ replaced by~$y$, and (for $i=1,2$)
\begin{eqnarray*}
&a_i(\cdot,\a):=a(x_i,\cdot,\a),\qquad
f_i(\cdot,\a):=f(x_i,\cdot,\a),&\\
&\ell_i(\cdot,\a):=\ell(x_i,\cdot,\a)+\tr(a(x_i,\cdot,\a)X_i)+f(x_i,\cdot,\a)\cdot p_i.&
\end{eqnarray*}
Taking into account the bounds
\begin{eqnarray*}
&\max|a_1-a_2|\leq L|x_1-x_2|,\qquad \max|f_1-f_2|\leq L|x_1-x_2|, &\\
&\max|\ell_1-\ell_2|\leq K(|X_1-X_2|+|p_1-p_2|)+L(1+|p_1|\wedge|p_2|+|X_1|\wedge|X_2|)|x_1-x_2|,&
\end{eqnarray*}
by Theorem~\ref{thm:cde}, we get
\begin{multline*}
\|v(\cdot;x_1,p_1,X_1)-v(\cdot;x_2,p_2,X_2)\|_\infty\leq
2C C_\ell L|x_1-x_2| +CK(|X_1-X_2|+|p_1-p_2|)\\
+CL(1+|p_1|\wedge|p_2|+|X_1|\wedge|X_2|)|x_1-x_2|.
\end{multline*}
Taking into account the inequality $C_\ell\leq (K+L+1) [1+(|p_1|+|X_1|)\wedge(|p_2|+|X_2|)]$,
 we accomplish the proof.
\end{Proofc}

Let us come back to the proof of Theorem~\ref{thm:hom}.
Fix $\e\in(0,1)$ and, for each~$\g\in(0,1)$, introduce the function
\begin{equation}\label{def:phi-x}
\phi(x):=u^\e(x)-u(x)-\e^2 v \left(\frac x\e;[u](x)\right)-\frac\g2|x|^2
\end{equation}
where $v(y;[u](x)):=v(y;x,Du(x),D^2u(x))$. Taking into account the bounds in Lemma~\ref{cor:p11}-$(i)$ and in Lemma~\ref{l3}, the functions $u^\e$, $u$ and $v(\cdot/\e;[u](\cdot))$ are bounded; whence  there exists a point~$\hx$ 
where the function $\phi$ attains its  maximum.

Set $c:=4C_2(1+2N)\e^\th$ (where $C_2$, $N$ and $\th$ are the constants introduced respectively in Lemma~\ref{cor:p11} and in Lemma~\ref{l3}) and introduce the function
\begin{equation}\label{def:B;tilde-phi-x}
\tilde\phi(x):=u^\e(x)-u(x)-\e^2 v \left( \frac x\e;[u](\hx)\right)-\frac\g2|x|^2-c|x-\hx|^2.
\end{equation}
The function~$\tilde\phi$ verifies $\tilde\phi(\hx)=\phi(\hx)$ and also, by the definition of $\hx$,
\begin{eqnarray*}\tilde\phi(\hx)-\tilde\phi(x)&=&[\phi(\hx)-\phi(x)]+[\phi(x)-\tilde\phi(x)]\geq
\phi(x)-\tilde\phi(x)\\
&\geq&-\e^2\left[v(x/\e;[u](x))-v(x/\e;[u](\hx))\right]+c\e^2
\end{eqnarray*}
for every $x\in \p B(\hx,\e)$.
Furthermore, owing to Lemma~\ref{cor:p11}-$(ii)$ and Lemma~\ref{l3}, it follows
\begin{equation*}
\tilde\phi(\hx)-\tilde\phi(x) \geq  -C_2[2N\e^\th +(1+\|Du\|_\infty+\|D^2u\|_\infty)\e]\e^2
+4C_2(1+2N)\e^{2+\th}> 0
\end{equation*}
for every $x\in \p B(\hx,\e)$.
Therefore, $\tilde \phi$ attains a maximum in some point $\tx\in B(\hx,\e)$.
Let us prove that there exists a constant $M_1>0$ (independent of $\e$ and $\g$) such that
\begin{equation}\label{cl:st-hx}
\g^{1/2}|\tx| \leq M_1.
\end{equation}
Actually, owing to Lemma~\ref{l3} and Lemma~\ref{cor:p11}-$(i)$, the inequality $\phi(\hx)\geq \phi(0)$ gives
\begin{equation*}
\frac \g2|\hx|^2 \leq 4N +2C_2(1+2N)\e^2.
\end{equation*}
For $M_1$ sufficiently large, we obtain: $\g^{1/2}|\hx| \leq M_1/2$; in particular, we deduce
$$
\g|\tx|\leq \g|\hx|+\g|\tx-\hx|\leq \g^{1/2}M_1/2+\g\e\leq \g^{1/2}M_1,
$$
which is our claim~\eqref{cl:st-hx}.

We claim now that there exists a constant~$M_2$ (independent of $\e$ and $\g$) such that
\begin{equation}\label{est:main-cl}
u^\e(\tx)-u(\hx)\leq M_2\left[\e^\th +\g^{1/2}\right].
\end{equation}
In order to prove this bound, we recall that the function~$u^\e$ solves problem~\eqref{homog} and that 
\[x\mapsto u^\e(x)-[u(x)+\e^2 v(x;[u](\hx))+\g|x|^2/2+c|x-\hx|^2]\]
attains a maximum in~$\tilde x$.
Thus, we have
\begin{multline*}
0\geq u^\e(\tx)+H\left(\tx,\tx/\e,Du(\tx)+\e D_y v(\tx/\e;[u](\hx))+\g\tx
+2c(\tx-\hx), \right.\\
\left. D^2u(\tx)+D^2_y v(\tx/\e;[u](\hx))+(\g+2c)I\right).
\end{multline*}
Assumptions~$(H1)$-$(H3)$, relation~\eqref{cl:st-hx}, Lemma~\ref{l3} and Lemma~\ref{cor:p11}-$(i)$ guarantee
\begin{align*}
& H\left(\tx,\tx/\e,Du(\tx)+\e D_y v(\tx/\e;[u](\hx))+\g\tx+2c(\tx-\hx),
 D^2u(\tx)+
 \right.\\ &\hskip 50mm\left.
+D^2_y v(\tx/\e;[u](\hx))+(\g+2c)I\right)\\
&\hskip 10mm \geq H\left(\tx,\tx/\e,Du(\tx),D^2u(\tx)
+D^2_y v(\tx/\e;[u](\hx))\right)-\\&\hskip 50mm-K\left[\e C_2(1+2N)+\g^{1/2}M_1+2c\e +\g+2c\right]\\
&\hskip 10mm \geq H\left(\hx,\tx/\e,Du(\hx),D^2u(\hx) +D^2_y
v(\tx/\e;[u](\hx))\right)-M_2\left[\e^\th +\g^{1/2}\right]
\end{align*}
for some constant~$M_2$ independent of $\e$ and $\g$.
Moreover, the cell problem~\eqref{s1p10} and the effective one~\eqref{eff} entail
\begin{align*}
&H\left(\hx,\tx/\e,Du(\hx),D^2u(\hx)+D^2_y v(\tx/\e;[u](\hx))\right)
 = \Ho(\hx,Du(\hx),D^2u(\hx))= -u(\hx).
\end{align*}
Substituting the last two relations in the previous one, we accomplish the proof of our claim~\eqref{est:main-cl}.

Finally, for every $x\in\R^n$, the inequalities $\tilde \phi(\tx)\geq \tilde\phi(\hx)=\phi(\hx)\geq \phi(x)$ give
\begin{eqnarray*}
u^\e(x)-u(x)&\leq& [u^\e(\tx)-u(\hx)]+[u(\hx)-u(\tx)]+\\
&&\hskip 20mm\e^2\left[v(x/\e;[u](x))-v(\tx/\e;[u](\hx))\right]+\frac \g2 |x|^2\\
&\leq& M_2\left[\e^\th +\g^{1/2}\right]+N\e +2C_2(1+2N)\e^2+\frac \g 2|x|^2
\end{eqnarray*}
where the latter inequality is due to relations~\eqref{est:main-cl}, to Lemma~\ref{l3} and to Lemma~\ref{cor:p11}-$(i)$.
Letting $\g\to 0^+$ (and increasing  $M_2$ if necessary), we deduce
\begin{equation*}
u^\e(x)-u(x)\leq M_2\e^\th\qquad \forall x\in\R^n.
\end{equation*}
Hence, we accomplished the proof of one inequality of the statement. Being similar, the proof of the other one is omitted.
\end{Proofc}
\begin{Remc}{Example}
Consider a diffusion coefficient~$a=a(x)$ with $a\in C^{1,1}(\R^n)$ satisfying $(H2)$-$(H3)$. Problem~\eqref{homog} is
\begin{equation}\label{ex}
u^\e-\tr(a(x)D^2 u^\e)+F_1\left(x,\frac x \e, D u^\e\right)=0
\end{equation}
where $F_1(x,y,p):=\max\limits_{\a\in A}\{-f(x,y,\a)\cdot p-\ell(x,y,\a)\}$.
Invoking \cite[Corollary 3.2]{AB8}, we have that the effective equation can be written as
\[u-\tr(a(x)D^2 u)+\bar F_1(x, D u)=0\]
where $\bar F_1(x,p):=\int_{(0,1]^n}F_1(x,y,p)\,dy$.
The regularity theory for semilinear equation (see~\cite{LU}) ensures that the effective solution~$u$ belongs to $C^{2,1}(\R^n)$. Therefore, by Theorem~\ref{thm:hom}, there exists a positive constant~$M$ such that
\[\|u^\e-u\|_\infty\leq M\e\qquad \forall \e\in(0,1).\]
\end{Remc}

\subsection*{Acknowledgment} The author is grateful to Professor P. Cardaliaguet, who suggested Theorem~\ref{thm:cdecont}, and to Professor E.R. Jakobsen for fruitful discussions.
The author wishes to thank the anonymous referees for pointing out a lack in the first version of the paper and for many useful remarks and suggestions.


\begin{thebibliography}{99}


\bibitem{AP11}
Y. Achdou,  and S. Patrizi.
\newblock{Homogenization of first-order equations with $u/\e$-periodic Hamiltonian: rate of convergence as $\e\to 0$ and numerical methods.}
\newblock{Math. Models Methods Appl. Sci.} {\bf 21} (2011), 1317 -- 1353. 

\bibitem{ACSppt}
S.N. Armstrong, P. Cardaliaguet and P.E. Souganidis.
\newblock{Error estimates and convergence rates for the stochastic homogenization of Hamilton-Jacobi equations.}
\newblock{Preprint available at http://arxiv.org/abs/1206.2601}.

\bibitem{AA}
V.I. Arnold and A. Avez.
\newblock{\em{Probl\`emes ergodiques de la m\`ecanique classique}}.
\newblock Gauthiers-Villars, Paris, 1967.


\bibitem{AB1}
O. Alvarez and M. Bardi.
\newblock{Viscosity solutions methods for singular perturbations in deterministic and stochastic control.}
\newblock{ SIAM J. Control Optim.} {\bf 40} (2001), 1159 -- 1188.



\bibitem{AB8}
O. Alvarez and M. Bardi.
\newblock{Ergodicity, stabilization, and singular perturbations for Bellman-Isaacs equation.}
\newblock{ Mem. Amer. Math. Soc.} {\bf 204} (2010), n. 960.

\bibitem{AL}
M. Arisawa and P.L. Lions,
\newblock{On ergodic stochastic control.}
\newblock{Comm. Partial Differential Equations} {\bf 23} (1998), 2187 -- 2217.




\bibitem{BL99} M. Bardi and F. Da Lio.
\newblock{On the strong maximum principle for fully nonlinear degenerate elliptic equations.}
\newblock{Arch. Math.} {\bf 73}  (1999), 276 -- 285.



\bibitem{BLM} G. Barles, O. Ley and H. Mitake.
\newblock{Short Time Uniqueness Results for Solutions of Nonlocal and Non-monotone Geometric Equations.}
\newblock{arXiv:1005.5597}.

\bibitem{BJ} G. Barles and E.R. Jakobsen.
\newblock{Error bounds for monotone approximation schemes for parabolic Hamilton-Jacobi-Bellman equations.}
\newblock{Math. Comp.} {\bf 76} (2007) 1861 -- 1893.

\bibitem{BS01} G. Barles and P. Souganidis.
\newblock{Space-time periodic solutions and long-time behavior of solutions to quasi-linear parabolic equations.}
\newblock{Siam J. Math. Anal.} {\bf 32} (2001), 1311--1323.



\bibitem{BLP}
A. Bensoussan, J.L. Lions and G. Papanicolaou.
\newblock {\em{Asymptotic Analysis for periodic Structures}}.
\newblock  {North-Holland, Amsterdam, 1978}.


\bibitem{BP99}
A. Bourgeat and A. Piatnitski.
\newblock{Estimates in probability of the residual between the random and the homogenized solutions of one-dimensional second-order operator.}
\newblock{ Asymptot. Anal.} {\bf 21} (1999), 303--315.



\bibitem{CS10}
L. Caffarelli and P. Souganidis.
\newblock {Rates of convergence for the homogenization of fully nonlinear uniformly elliptic pde in random media.}
\newblock{Invent. Math.} {\bf 180} (2010), 301 -- 360.



\bibitem{CSW}
L. Caffarelli, P. Souganidis and L. Wang.
\newblock Homogenization of fully nonlinear, uniformly elliptic and parabolic partial differential equations in stationary ergodic media.
\newblock { Comm. Pure Appl. Math.} {\bf 58} (2005), 319 -- 361.

\bibitem{CCM11}
F. Camilli, A. Cesaroni and C. Marchi.
\newblock{Homogenization and vanishing viscosity in fully nonlinear elliptic equations: rate of convergence estimates.}
\newblock Adv. Nonlinear Stud. {\bf 11} (2011), 405 -- 428.

\bibitem{CM09}
F. Camilli and C. Marchi.
\newblock{Rates of convergence in periodic homogenization of fully nonlinear uniformly elliptic PDEs.}
\newblock{Nonlinearity} {\bf 22} (2009), 1481 -- 1498.

\bibitem{CM11}
F. Camilli and C. Marchi.
\newblock{On the convergence rate in multiscale homogenization of fully nonlinear elliptic problems.}
\newblock{Netw. Heterog. Media} {\bf 6} (2011), 61 -- 75.

\bibitem{CI01}
I. Capuzzo Dolcetta and H. Ishii.
\newblock{On the rate of convergence in homogenization of Hamilton-Jacobi equations.}
\newblock Indiana Univ. Math. J. {\bf 50} (2001), 1113 -- 1129.


\bibitem{CFS}
I.P. Cornfeld, S.V. Fomin and Y.G. Sinai.
\newblock {\em{Ergodic theory}}.
\newblock  {Springer-Verlag, Berlin, 1982}.





\bibitem{DK}H. Dong and N.V. Krylov.
\newblock{The rate of convergence of finite-difference approximations for parabolic Bellman equations with Lipschitz coefficients in cylindrical domains.}
\newblock{Appl. Math. Optim.} {\bf 56} (2007), 37 -- 66.



\bibitem{Ev2}
L. Evans.
\newblock{Periodic homogenisation of certain fully nonlinear partial differential equations}
\newblock{Proc. Roy. Soc. Edinburgh Sect. A} {\bf 120} (1992), 245 -- 265.


\bibitem{FlS}
W.H. Fleming and H.M. Soner.
\newblock{\em Controlled Markov processes and viscosity solutions.}
\newblock {Springer-Verlag}, Berlin 1993.



\bibitem{GT}
D. Gilbarg and N.S. Trudinger
\newblock{\em Elliptic Partial Differential Equations of Second Order.}
\newblock{Springer}, Berlin 1983. Second edition.





\bibitem{JG}
E.R. Jakobsen and C.A. Georgelin.
\newblock{Continuous dependence results for non-linear Neumann type boundary value problems.}
\newblock{J. Differential Equations} {\bf 245} (2008), 2368 -- 2396.

\bibitem{JK1}
E.R. Jakobsen and K.H. Karlsen.
\newblock{Continuous dependence estimates for viscosity solutions of fully nonlinear degenerate parabolic equations.}
\newblock{J. Differential Equations} {\bf 183} (2002), 497 -- 525.



\bibitem{JKO}
V.V. Jikov, S.M. Kozlov and O.A. Oleinik.
\newblock{\em{Homogenization of Differential Operators and Integral Functionals}}.
Springer, Berlin, 1994.






\bibitem{LU} O.A. Ladyzenskaja and N.N. Ural'tseva.
\newblock Linear and quasilinear elliptic equations.
\newblock {\em{Linear and quasilinear elliptic equations}}.
Academic Press, New York, 1968.






\bibitem{M12}
C. Marchi.
\newblock{Continuous dependence estimates for the ergodic problem of Bellman-Isaacs operators via the parabolic Cauchy problem.}
\newblock{ESAIM Control Optim. Calc. Var.} {\bf 18} (2012), 954 -- 968.

\bibitem{Sa88}
M.V. Safonov.
\newblock{On the classical solution of nonlinear elliptic equations of second-order.}
\newblock{Math. USSR-Izv.} {\bf 33} (1989), 597 -- 612. (Engl. transl. of \newblock{Izv. Akad. Nauk SSSR Ser.
Mat.} {\bf 52} (1988)).





\bibitem{t1}
N.S. Trudinger.
\newblock{Comparison principles and pointwise estimates for viscosity solutions of nonlinear elliptic equations.}
\newblock{Rev. Mat. Iberoamericana} {\bf 4} (1988), 453 -- 468.

\bibitem{t2}
N.S. Trudinger.
\newblock{On regularity and existence of viscosity solutions of nonlinear second order, elliptic equations.}
\newblock{{\em Partial differential equations and the calculus of variations} (Eds.: F. Colombini, A. Marino, L. Modica and S. Spagnolo).}
Birkh\"auser, Boston, 1989.



\bibitem{Y86}
V.V. Yurinski.
\newblock{Averaging of symmetric diffusion in a random media.}
\newblock Sibirsk. Mat. Zh. {\bf 27} (1986), 167 -- 180; Engl. transl.
Siberian Math. J. {\bf 27} (1986), 603 -- 613.


\end{thebibliography}
\end{document}